\documentclass[12pt]{amsart}
\usepackage{amsmath,amsthm}

\def\pmatrix{\left(\begin{matrix}}
\def\endpmatrix{\end{matrix}\right)}

\def\Sp{\operatorname{Sp}(g,\Z)}

\def\Z{{\mathbb Z}}
\def\F{{\mathbb F}_2}
\def\C{{\mathbb C}}

\def\PP{{\mathbb P}}

\def\de{\delta}

\def\p{\partial}

\def\t{\theta}

\def\T{\Theta}
\def\e{\varepsilon}

\def\o{\overline}
\def\A{{\mathcal A}}
\def\M{{\mathcal M}}

\def\H{{\mathcal H}}
\def\T{{\mathcal T}}

\def\tch#1#2{{\left[\begin{matrix}#1\\ #2\end{matrix}\right]}}
\def\tt#1#2{{\t\tch{#1}{#2}}}

\theoremstyle{plain}
\newtheorem{thm}{Theorem}%[section]
\newtheorem{lm}[thm]{Lemma}
\newtheorem{prop}[thm]{Proposition}
\newtheorem{cor}[thm]{Corollary}

\theoremstyle{definition}

\newtheorem{rem}[thm]{Remark}
\begin{document}
\title[The cosmological constant and the Schottky form]{The superstring cosmological constant and the Schottky form in genus 5}
\author{Samuel Grushevsky}
\address{Mathematics Department, Princeton University, Fine Hall,
Washington Road, Princeton, NJ 08544, USA. }
\thanks{Research of the first-named author is supported in part by National Science Foundation under the grant DMS-05-55867.}
\email{sam@math.princeton.edu}
\author{Riccardo Salvati Manni}
\address{Dipartimento di Matematica, Universit\`a ``La Sapienza'',
Piazzale A. Moro 2, Roma, I 00185, Italy}
\email{salvati@mat.uniroma1.it}
\date{\today}
\begin{abstract}
Combining certain identities for modular forms due to Igusa with Schottky-Jung relations, we study the cosmological constant for the recently proposed ansatz for the chiral superstring measure in genus 5. The vanishing of this cosmological constant turns out to be equivalent to the long-conjectured vanishing of a certain explicit modular form of genus 5 on the moduli of curves $\M_5$, and we disprove this conjecture, thus showing that the cosmological constant for the proposed ansatz does not vanish identically. We exhibit an easy modification of the genus 5 ansatz satisfying factorization constraints and yielding a vanishing cosmological constant. We also give an expression for the cosmological constant for the proposed ansatz that should hold for any genus if certain generalized Schottky-Jung identities hold.
\end{abstract}
\maketitle
\section{Introduction}
\subsection*{Superstring scattering amplitudes}

In (super)string perturbation theory one of the main questions is finding an explicit expression for the string measure to an arbitrary loop (genus) order. For the bosonic string, the infinite-dimensional integral over the space of possible worldsheets becomes, after taking a quotient by the appropriate group of symmetries, an integral over the moduli space of curves $\M_g$. For superstrings, computing the measure is much more difficult, as the integration over the odd variables needs to be performed. However, for $g=1$ this difficulty is not present, and an expression for the genus 1 superstring measure was obtained by Green and Schwarz \cite{grsc}. In general the superstring measure is a measure on the moduli space of super-Riemann surfaces, and dealing with the supermoduli and finding an appropriate gauge-fixing procedure is extremely difficult. In a series of papers \cite{DHP1,DHPa,DHPb,DHPc,DHPd,DHPe} D'Hoker and Phong overcame these conceptual difficulties, derived an expression for the chiral superstring measure for $g=2$ from first principles, and verified its further physical properties.

In \cite{159,182} D'Hoker and Phong establish the modern program of finding higher genus (i.e. multiloop) superstring amplitudes by using factorization constraints: they investigated the behavior of the amplitude as the curve degenerates, and showed that in the  limit the constant term is expected to factorize as the product of lower-dimensional amplitudes. D'Hoker and Phong then proposed and studied an expression for the genus 3 superstring amplitude assuming holomorphicity of certain square roots. Cacciatori, Dalla Piazza, and van Geemen in \cite{CDPvG} advanced the program begun by D'Hoker and Phong by constructing a holomorphic modular form in genus 3 satisfying the factorization constraints, while in \cite{DPvG} the uniqueness of this ansatz is shown.

\smallskip
The problem of finding the superstring measure for higher genera presents additional difficulties. While for $g\le 3$ loops the moduli space of curves $\M_g$ is dense in the moduli space $\A_g$ of ppavs, this is no longer the case for genus $g\ge 4$. Thus while for $g\le 3$ an ansatz for a superstring measure is a Siegel modular form on $\A_g$, for $g\ge 4$ such an ansatz for the amplitude may only be defined on $\M_g$, and not on all of $\A_g$.

In \cite{GR} the first-named author presented a general framework for the ansatz\"e of Green-Schwartz, D'Hoker-Phong, and Cacciatori-Dalla Piazza-van Geemen, proposed an ansatz for genus 4, and a possible ansatz, satisfying the factorization constraints, for the superstring measure in any genus, subject to the condition that certain holomorphic roots of modular forms are single-valued on $\M_g$ (the genus 4 ansatz was then also obtained independently in \cite{CDPvG2}). In \cite{RSM} the second-named author proved that in genus 5 these holomorphic square roots are indeed well-defined modular forms on a suitable covering of $\M_5$, and thus that an ansatz is well-defined for $g=5$. For a review and further developments, see Morozov \cite{M1,M2}. See also \cite{MV1} for a different approach.

In \cite{GSM} we showed that in genus 3 the 2-point function vanishes as expected. Recently in \cite{MV} Matone and Volpato showed that certain quantities connected with this ansatz and the 3-point function no longer vanish for genus 3. However, they discuss also the possibility of a non-trivial correction which would still result in the vanishing of the 3-point function, and further investigation of the question is required. Whatever correction may be required will not influence, however, the first constraint on the ansatz --- the vanishing of the cosmological constant along $\M_g$. This has been verified in \cite{GR} and \cite{RSM} for all $g\le 4$ by using an expression for the ansatz in terms of theta series associated to quadratic forms.

\subsection*{The Schottky problem}

In complex algebraic geometry, the Riemann-Schottky problem --- the question of determining which ppavs are Jacobians of curves --- is one of the classical problems in algebraic geometry. The solution to the problem for the first non-trivial case, that of $g=4$, was given by Schottky in \cite{sc}, who constructed an explicit modular form vanishing on the locus of Jacobians of curves of genus 4. Igusa in \cite{Ch} proved that the Schottky form is irreducible and thus defines $\M_4\subset\A_4$.

The original expression for the Schottky form in \cite{sc} used the combinatorics of theta characteristics, but it was shown by Igusa that it is equal to the weight 8 modular form $F_g:=(2^g-1)\sum\theta_m^{16}-2(\sum\theta_m^8)^2$, for genus $g=4$
(where the sum is taken over all even theta characteristics $m$). It is known that for $g\le 3$ the modular form $F_g$ is identically zero on $\A_g$ --- this follows from Riemann's bilinear relations for $g=1,2$, and is the only equation relating theta constants of the second order for $g=3$, cf. \cite{vgvdg}. The properties of the modular form $F_g$ for higher genus are of obvious interest, and some identities for it were obtained by Igusa in \cite{Ch}.

In the 1980s Belavin, Knizhnik, and Morozov \cite{BK,BK1,M86}, and later D'Hoker and Phong \cite{159} for physics reasons conjectured that the form $F_g$ vanishes on $\M_g$ for any $g$ (it was shown by Igusa that $F_g$ is proportional to $f_4^2-f_8$, where $f_4$ and $f_8$ are the theta series associated to the even lattices $E_8\times E_8$ and $D_{16}^+$ --- see below for the definitions --- and physically the vanishing of $F_g$, i.e. the equality of $f_4^2$ and $f_8$, would be interpreted as the equality of the measures for the $SO(32)$ and $E_8\times E_8$ theories). This question was investigated by Poor who in \cite{Po} showed that for any $g$ the form $F_g$ vanishes on the locus of hyperelliptic curves; the conjectural vanishing of $F_g$ on $\M_g$ remained open.

\subsection*{Summary of results}

In this paper we study the cosmological constant for the proposed in \cite{GR} chiral superstring scattering measure for $g\ge 5$. By extending the results of Igusa \cite{Ch} on certain modular forms constructed using theta constants associated to syzygetic subspaces of characteristics we obtain an alternative expression for one of the summands in the ansatz. This allows us to show that in genus 5 the cosmological constant for the proposed ansatz is in fact equal to a non-zero multiple of the Schottky form $F_5$ described above. Thus the vanishing of the cosmological constant for the proposed genus 5 ansatz is equivalent to the  conjecture of Belavin, Kniznik, Morozov, D'Hoker, and Phong that $F_5$ vanishes identically on $\M_5$.

By studying the boundary of $\M_5$ and computing the degenerations of the proposed ansatz, using theta functional techniques, we show that this conjecture is {\it false}: that $F_5$ in fact does {\it not} vanish identically on $\M_5$. By using the results on the slope of effective divisors on $\M_g$ we identify the zero locus of $F_5$ in $\M_5$ as the divisor of trigonal curves in $\M_5$.

Thus it follows that the cosmological constant for the originally proposed in \cite{GR} ansatz for the chiral superstring measure does not vanish identically in genus 5. We thus give a simple explicit formula for a modification of the genus 5 ansatz that satisfies factorization constraints and results in an identically vanishing cosmological constant.

It follows from our results that $F_g$ does not vanish identically on $\M_g$ for any $g\ge 5$. Subject to the validity of certain generalized Schottky-Jung identities involving roots of degree $2^{g-4}$ of polynomials in theta constants (see the appendix), we further show that the cosmological constant for the originally proposed ansatz is equal to a non-zero multiple of $F_g$  in any genus, and thus also does not vanish identically on $\M_g$ for any $g\ge 5$ --- thus an adjustment of the ansatz and a further investigation of the question would be of interest.

\section*{Acknowledgements}
We would like to thank Sergio L. Cacciatori, Francesco Dalla Piazza, Bert van Geemen, Marco Matone, Duong Phong, and Roberto Volpato for useful discussions and suggestions regarding chiral superstring scattering amplitudes.

We would like to thank Giulio Codogni for detailed analysis of our argument and for pointing out the mistake in our original proof of theorem 15, and we are grateful to Enrico Arbarello for stimulating discussions.

\section{Notations}
We denote by $\H_g$ the Siegel upper half-space of symmetric complex matrices with positive-definite imaginary part, called period matrices. The action of the symplectic group $\Sp$ on $\H_g$ is given by
$$
 \pmatrix A&B\\ C&D\endpmatrix\circ\tau:= (A\tau+B)(C\tau+D)^{-1}
$$
where we think of elements of $\Sp$ as of consisting of four
$g\times g$ blocks, and they preserve the symplectic form given in
the block form as $\pmatrix 0& 1\\ -1& 0\endpmatrix$.

For a period matrix $\tau\in\H_g$, $z\in \C^g$ and $\e,\de\in \F^g$ (where $\F$ denotes the abelian group $\Z/2\Z=\lbrace 0,1\rbrace$ for which we use the additive notation)
the associated theta function with characteristic $m=[\e, \de]$ is
$$
 \t_m(\tau, z)=\tt\e\de(\tau,z)=\sum\limits_{n\in\Z^g}\exp(\pi i ((n+\e/2)'\tau
 (n+\e/2)+ 2(n+\e/2)'( z+\de/2))
$$
(where we denote by $X'$ the transpose of $X$). As a function of $z$, $\t_m(\tau, z)$ is odd or even depending on whether
the scalar product $\e\cdot\de\in\F$ is equal to 1 or 0,
respectively. Theta constants are restrictions of theta
functions to $z=0$. We shall write $\t_m$ for theta constants.\smallskip

For a set of characteristics $M=(m_1, m_2,\dots, m_k)$ we set
$$
 P(M):=\prod_{i=1}^k\t_{m_i}.
$$
A holomorphic function $f:\H_g\to\C$ is a modular form
of weight $k/2$ with respect to a subgroup $\Gamma\subset\Sp$ of finite index if
$$
 f(\gamma\circ\tau)=)\det(C\tau+D)^{k/2}f(\tau)\quad
 \forall\gamma\in\Gamma,\forall\tau\in\H_g.
$$
and if additionally $f$ is holomorphic at all cusps when $g=1$.
We denote by $[\Gamma, k/2]$ the vector space of such functions.
Theta constants are modular forms of weight $k/2$ with respect to a certain subgroup $\Gamma(4,8)$. For further use, we denote $\Gamma_g:=\Sp$ the integral symplectic group and $\Gamma_g (1,2)$ its subgroup defined by
$${\rm diag}(AB')\equiv {\rm diag}(CD')\equiv 0\,\, mod 2.$$

\section{The ansatz for the chiral superstring measure}
We recall from \cite{GR} or \cite{RSM}
\begin{lm}
If $16$ divides $2^i s$, and $g\geq i$, then
$$
 P_{i, s}^g (\tau):= \sum_V P(V)^s(\tau),
$$
where the sum is over all $i$-dimensional subspaces of $\F^{2g}$,
belongs to $[\Gamma_g(1,2),\, 2^{i-1}s]$.
\end{lm}

For any even characteristic $m$, we can define $P_{i, s}^g[m] (\tau)$ by taking the sum above over all affine subspaces $V$ of $\F^{2g}$ of dimension $g$ (i.e. translates of $i$-dimensional linear subspaces) containing the characteristic $m$. The function $P_{i, s}^g[m] (\tau)$ is then a modular form with respect to a subgroup of $\Gamma_g$ conjugate to $\Gamma_g (1,2)$ (note that $\Gamma_g (1,2)\subset\Gamma_g$ is not normal).
\begin{cor}
If $16$ divides $2^i s$, and $g\geq i$, the form
$$
 S_{i,s}^g:=\sum_{m\in\F^{2g}}\sum_V P(V+m)^s(\tau)=\sum_m P_{i,s}^g[m](\tau),
$$
where the sum is taken over all $i$-dimensional linear subspaces $V$, belongs to $[\Gamma_g,\, 2^{i-1} s]$.
\end{cor}

We recall the main results obtained in \cite{GR} and \cite{RSM}.
\begin{prop}[\cite{GR}]
The modular forms $P_{i, s}^g$ restrict to the locus of block diagonal period matrices $\H_k\times\H_{g-k}$ as follows:
$$
P_{i, s}^g\left(\begin{matrix}\tau_1&0\\ 0&\tau_2\end{matrix}\right) = \sum_{0\leq n,m\leq i\leq n+m}N_{n,m;i} P_{n, 2^{i-n}s}^k(\tau_1) P_{m, 2^{i-m}s}^{g-k}(\tau_2),
$$
where
$$N_{n,m;i}:=\prod_{j=0}^{n+m-i-1}\frac{(2^n -2^j)(2^m-2^j)}{2^{n+m-i }-2^j},$$
for any $\tau_1\in\H_k$ and $\tau_2\in\H_{g-k}$.
\end{prop}
\begin{thm}[\cite{GR}]
For $g \leq 4$ the function
$$
 \Xi^{(g)}[0]:=\frac{1}{2^g}\sum\limits_{i=0}^g (-1)^i2^{\frac{i(i-1)}{2}}P_{i, 2^{4-i}}^{g}
$$
is a modular form in $[\Gamma_g(1,2),8]$, and its restriction to $\H_k\times\H_{g-k}$ is
$$
 \Xi^{(g)}[0]\left(\begin{matrix}\tau_1&0\\ 0&\tau_2\end{matrix}\right)=\Xi^{(k)}[0](\tau_1)\cdot\Xi^{(g-k)}[0](\tau_2),
$$
for any $\tau_1\in\H_k,\tau_2\in\H_{g-k}$.
\end{thm}
We also define for any characteristic $m$
$$
 \Xi^{(g)}[m]:=\frac{1}{2^g}\sum\limits_{i=0}^g (-1)^i 2^{\frac{i(i-1)}{2}}P_{i, 2^{4-i}}^{g}[m].
$$
These $\Xi^{(g)}[m]$ satisfy similar factorization constraints, and thus are natural candidates for the chiral superstring measure.

In \cite{GR} it is also shown that the above statement holds for $g > 4$ as well, up to a possible inconsistency in the modularity. In fact, since $\H_g$ is simply connected, the individual degree $2^n$ roots needed to define the $P$'s above are well-defined globally, but they are not necessarily modular forms due to possible sign inconsistency.
\smallskip

We denote $\A_g:=\H_g/\Gamma_g$ the moduli space of principally polarized abelian varieties. The Torelli map gives an immersion of the moduli space of curves $\M_g\hookrightarrow\A_g$. We define $\T_g\subset\H_g$ (the Torelli space) as the preimage of $\M_g\subset\A_g$ under the projection $\H_g\to\A_g$.

\begin{thm}[\cite{RSM}]
The restrictions to $\T_5$ of $P_{5, 1/2}^{5}$ (which are sums of square roots of polynomials in theta constants) satisfy the modularity condition with respect to $\Gamma_5(1, 2)$. Hence $\Xi^{(5)}[0]|_{\T_5}$ is a section of the restriction of the bundle of modular forms of weight 8 to $\T_5$.
\end{thm}

Summing up we get the following expression for the cosmological constant
\begin{equation}\label{Xi}
 \Xi^{(g)}:=\sum_m \Xi^{(g)}[m] =\frac{1}{2^g}\sum\limits_{i=0}^g (-1)^i 2^{\frac{i(i+1)}{2}}S_{i, 2^{4-i}}^g
\end{equation}
(note that when summing over characteristics $m$, each polynomial $P_i$ appears in $2^i$ terms; thus we have the $2^{\frac{i(i+1)}{2}}$ in the formula for $\Xi^{(g)}$ here instead of the $2^{\frac{i(i-1)}{2}}$ in the formula for $\Xi^{(g)}[m]$).
As a consequence of the previous discussion we have
\begin{prop}
When $g\leq 4$, $\Xi^{(g)}$ is a modular form on $ \H_g$ of weight 8, while
$\Xi^{(5)}|_{\T_5}$ is a section of the restriction of the line bundle of modular forms of weight 8.
\end{prop}

For any even positive definite unimodular matrix $S$ of degree $2k$ (i.e $S$ is positive definite, $\det S=1$ $x'Sx\equiv 0 \forall x \in \Z^k$), we define the theta series for $\tau\in\H_g$ by
$$
 f_S^{(g)}(\tau):=\sum\limits_{u\in\Z^{k,g}}\exp(\pi i {\rm tr}(u'Su\tau))
$$
These are modular forms in $[\Gamma_g,k]$, cf.~\cite{Fr}.

For any $g$, we denote by $f_4^{(g)}$ and $f_8^{(g)}$ the theta series, of weights 4 and 8, respectively, associated to the even unimodular matrices related to the lattices $E_8$ and $D_{16}^+$, cf.~\cite{CS}. It is shown in \cite{Ch} that the following identity holds for any $g$:
\begin{equation}\label{fS}
 ((f_4^{(g)})^2-f_8^{(g)})={\frac{1}{2^{2g}}}((1-2^g )S_{0,16}^g+2S_{1,8}^g).
\end{equation}
\smallskip
In this note we consider the case $g=5$ of the ansatz for chiral superstring measures, which involves square roots of theta constants.
In this case we prove
$$
 \Xi^{(5)}= \frac{-51}{217}((f_4^{(5)})^2-f_8^{(5)})
$$
on $\T_5$, and then resolve a long-standing open question posed in \cite{BK,BK1,M86,159}, by showing that this expression does {\it not} vanish identically. In an appendix we then discuss a possible generalization of the Schottky-Jung identities, and conjectural results for arbitrary genus.

\section{ The cosmological constant }
The vanishing of the cosmological constant means for $\Xi^{(g)}(\tau)$ to vanish identically on the moduli space of curves $\M_g$ (we will work on $\T_g$). This has been verified for the proposed ansatz in genus 2 in \cite{DHP1}, for genus 3 in \cite{CDPvG}, and for genus 4 in \cite{GR} and \cite{RSM}, so that we know that $\Xi^{(g)}(\tau)$ is identically zero for $\tau\in\T_g$ for $g\leq 4$. The proof given in \cite{RSM} was an immediate consequence of remarkable formulas deduced from the Riemann relations by Igusa around thirty years ago.
\begin{lm}[\cite{Ch}]\label{iglm}
We have
$$(2^{2g}-1)S_{0,16}^g= 6S_{1,8}^g+24S_{2,4}^g$$
$$(2^{2g-2}-1)S_{1,8}^g= 18S_{2,4}^g+168S_{3,2}^g$$
$$(2^{2g-4}-1)S_{2,4}^g= 42S_{3,2}^g+840 S_{4,1}^g$$
for $g\geq 2,3,4$ respectively.
\end{lm}
In this paper we study the cosmological constant when $g\ge 5$, and to this end recall the form of the Riemann relations used to prove the above result. For any $ a\in \F^{2g}$ we denote by $a'$ (resp. $a''$) the vectors of the first $g$ (resp. last $g$) entries. For any $a, b, c\in \F^{2g}$ we set
$$
 (a,b,c)=\exp(\pi i\sum_{j=1}^g (a'_j b''_j c''_j +a''_j b'_j c''_j+ a''_j b''_j c'_j )
$$
and observe that this is a symmetric tricharacter.\smallskip

We also set $e(a,b):=(a,a,b)(a,b,b)$; with these notations, Riemann's theta formula can be stated as follows

\begin{lm}[\cite{Ch}]
For any $m, a, b$ in $ \F^{2g}$ we have
\begin{equation}\label{riemannrelation}
 (m,a,a)(m,b,b)(m,a,b)\t_m\t_{m+a}\t_{m+b}\t_{m+a+b}
\end{equation}
$$
 =2^{-g}\sum_{n\in\F^{2g}}e(m,n)(n,a,a)(n,b,b)(n,a,b)\t_n\t_{n+a}\t_{n+b}\t_{n+a+b}.
$$
\end{lm}
We remark that our proof of the modularity of $\Xi^{(5)}$ also uses the so called Schottky-Jung relation for theta constants of Jacobians of curves, cf.~\cite{RF,vG, Ts}. It is well-known that Riemann relations in genus $g$ induce Schottky relations for periods of Jacobians in genus $g+1$ of similar structure: Riemann relations involve homogeneous monomials of degree 4 in the $\t_m$, while the Schottky-Jung relations involve monomials of degree 8 in the square root of $\t_m$. To write them down explicitly, we set $c:=\tch{0&\ldots&0&0}{0&\ldots&0&1}\in\F^{2g+2}$, and for any $m=\tch{m'}{m''}\in\F^{2g}$ denote $\o{m}:=\tch{m'&0}{m''&0}\in\F^{2g+2}$. Then as a consequence of  the classical Schottky-Jung relations we have
\begin{lm}
For any $\pi\in\T_{g+1}$, the following identity holds for theta constants evaluated at $\pi$
$$
(\o{m},\o{a},\o{a})(\o{m},\o{b},\o{b})(\o{m},\o{a},\o{b})\sqrt{\t_{\o{m}}\t_{\o{m+a}}\t_{\o{m+b}}\t_{\o{m+a+b}}
\t_{\o{m}+c}\t_{\o{m+a}+c}\t_{\o{m+b}+c}\t_{\o{m+a+b}+c}}
$$
$$
 =2^{-g}\sum_{n\in \F^{2g}}e(\o m, \o n)(\o n,\o a,\o a)(\o n,\o b,\o b)(\o n,\o a,\o b)$$
 $$\sqrt{\t_{\o{n}}\t_{\o{n+a}}\t_{\o{n+b}}\t_{\o{n+a+b}}
\t_{\o{n}+c}\t_{\o{n+a}+c}\t_{\o{n+b}+c}\t_{\o{n+a+b}+c}}
$$
\end{lm}
We observe that the structure of Riemann and Schottky-Jung relations  above is the same, and that Riemann relations in genus $g$ and Schottky relations in genus $g+1$ have the same coefficients. Thus the $g=4$ Riemann relation (valid on $\H_4$) of the form
$$
 r_1\pm r_2\pm r_3\pm r_4=0,
$$
with each $r_i$ a monomial of degree 4 in theta constants,
induces a Schottky relation in genus $5$, valid on $\T_5$, of the form
$$
 \sqrt R_1\pm \sqrt R_2\pm \sqrt R_3\pm \sqrt R_4=0,
$$
where each $R_i$ is the square root of a monomial of degree 8 in theta constants of a Jacobian in $\T_5$, cf.~\cite{AC}.

Note that if we replace $\o n$ with $\o n +c$ in the right-hand-side of the lemma above, nothing changes, since
$e(\o m, \o n +c)=e(\o m, \o n)$ and $(\o n,\o a,\o b)=(\o n +c,\o a,\o b)$. Moreover, if we let $d:=\tch{0&\ldots&0&1}{0&\ldots&0&0}\in\F^{2g+2}$,
then for any $n\in\F^{2g}$ one of the characteristics $\o{n}+d$ and $\o{n}+c+d$ is odd, and thus  $\t_{\o{n}+d}\t_{\o{n}+c+d}=0$. Thus in the lemma above we can extend the summation over all of $\F^{2g+2}$ to get

\begin{lm}
For any $\pi\in\T_{g+1}$, the following identity holds for theta constants evaluated at $\pi$
$$
(\o{m},\o{a},\o{a})(\o{m},\o{b},\o{b})(\o{m},\o{a},\o{b})\sqrt{\t_{\o{m}}\t_{\o{m+a}}\t_{\o{m+b}}\t_{\o{m+a+b}}
\t_{\o{m}+c}\t_{\o{m+a}+c}\t_{\o{m+b}+c}\t_{\o{m+a+b}+c}}
$$
$$
 =2^{-g-1}\sum_{n\in \F^{2g+2}}e(\o m, n)( n,\o a,\o a)( n,\o b,\o b)( n,\o a,\o b)$$
 $$\sqrt{\t_{n}\t_{n+\o a}\t_{n+\o{b}}\t_{n+\o{a+b}}
\t_{n+c}\t_{n+\o{a}+c}\t_{n+\o{b}+c}\t_{n+\o{a+b}+c}}
$$
\end{lm}

All the terms appearing in the above   relation   are of the form
$\sqrt{P(N+n)}$, with $N=\langle \o a, \o b , c\rangle$ (we denote by $\langle\ \rangle$ the linear span). Such a polynomial $P(N+n)$ is not identically zero if and only if $N+n$ is an even coset of a totally isotropic, with respect to the form  $e(m,n)$, 3-dimensional space $N$. The symplectic group acts transitively on such cosets, and $\sigma(N)P(N+n)(\tau)$ maps to  $\sigma(N_1)P(N_1+n_1)(\tau)$, where the sign $\sigma(N) =\pm 1$  depends only on the subspace $N$ and not on the coset, cf.~\cite{Ch}.

As an immediate consequence  we get a more general result than the above lemma that can be stated for any totally isotropic space $N$. Since square roots appear in the formula, and there is a choice of a sign for each of them, there will be signs $\sigma(n)$ depending on the cosets. Applying the same argument as in
\cite{RSM}, we get the following constraint on the signs:
$$
 \sigma(n_1)\sigma(n_2)\sigma(n_3)\sigma(n_4)=1
$$
if $n_1+n_2+ n_3+n_4 =0$.

\begin{lm}
Let $N=\langle a, b, c\rangle$ be a totally isotropic 3-dimensional subspace of $\F^{2g+2}$. Then for any $\pi\in\T_{g+1}$, the following identity holds for theta constants evaluated at $\pi$
$$
(m,a,a)(m,b,b)(m,a,b)\sigma(m)\sqrt{P(N+m)}=
$$
$$2^{-g-1}\sum_{n\in \F^{2g+2}}e( m, n)\sigma(n)( n,a,a)( n,b,b)( n,a,b)\sqrt{P(N+n)}
$$
\end{lm}
To obtain special relations for theta constants of Jacobians using the above Schottky-Jung identity, we proceed as in \cite{Ch}, where identities for theta constants of arbitrary abelian varieties are obtained by using Riemann relations. We repeat the argument given there, since it is elementary, but quite involved.

We take the fourth  powers of both sides of the above formula and sum over $m\in\F^{2g+2}$ to obtain
$$
 \sum_mP(N+m)^2=2^{-2g-2}\left(\sum_nP(N+n)^2+3!\sum_{n,m}P(N+n)P(N+m)\right.
$$
$$
  \left.+4!\sum_{n_1,n_2,n_3,n_4}\sqrt{P(N+n_1)P(N+n_2)P(N+n_3)P(N+n_4)}\right)
$$

Because of the orthogonality of the characters, the only non-zero terms here would be the ones with $n_1+n_2+ n_3+n_4 =0$. The crucial observation is that the signs $\sigma$ disappear in this formula, so that it now looks exactly similar to the one for theta constants of arbitrary abelian varieties, given in \cite{Ch}. Indeed, in the first term two terms on the right we have $\sigma^2=1$, while in the last term we have $\sigma(n_1)\sigma(n_2)\sigma(n_3)\sigma(n_4)$, which is equal to 1 since $n_1+n_2+n_3+n_4=0$.

\smallskip
We now sum over all 3-dimensional isotropic subspaces $N$, and note that
$$\sum_{N}\sum_{m} \prod _{n\in Nm}\t_n^2= \sum_{N}\sum_{m}P(Nm)^2=8 S_{3,2}^{g+1}$$
where the coefficient $8$ is due to the fact that each $N$ has 8 elements. We further compute
$$
 \sum_{N}\sum_{n,m }P (N+m)  P(N+n) = 28 S_{3,2}^{g+1} + 64\cdot 15 S_{4,1}^{g+1},
$$
where $28$ appears as the number of pairs of distinct elements of $N$, $64=8\cdot 8$ is the  number of elements in $(N+n)\times (N+m)$, and $15$ is the  number of 3-dimensional isotropic spaces contained in a fixed 4-dimensional isotropic space (the 4-dimensional space $(N+n)\sqcup(N+m)$ must be isotropic, otherwise the corresponding product is zero).

Finally for the last term we get
$$
\sum_{N} \sum_{ n_1, n_2, n_3, n_4 }\quad \sqrt{P(N+n_1)P(N+n_2)P(N+n_3)P(N+n_4)}=
$$
$$
14  S_{3,2}^{g+1} + 112\cdot 15 S_{4,1}^{g+1} +2^9\cdot 155 S_{5, 1/2}^{g+1}.
$$
Here $14$ is the  numbers of quadruplets $ n_1, n_2, n_3, n_4$ such that $n_1+n_2+ n_3+n_4 =0$, and all $n_i$ are in $N$, $112= 28\cdot 4$ is the  number of quadruplets $ n_1, n_2, n_3, n_4$ such that $n_1+n_2+ n_3+n_4 =0$ with $n_1, n_2\in N$ and $n_3, n_4\in N+n_3$,
$2^9$ is the  numbers of quadruplets $ n_1, n_2, n_3, n_4$ such that $n_1+n_2+ n_3+n_4 =0$ with all $N+n_i$ disjoint, and $155$ is the  number of 3-dimensional isotropic spaces contained in a 5-dimensional isotropic space (notice that in this case the union $\sqcup_i(N+n_i)$ is 5-dimensional isotropic).

\smallskip
Applying these results, we finally get the following:
$$
 8 S_{3,2}^{g+1}={2^{-2g-2}}\left(8 S_{3,2}^{g+1}+ 6( 28 S_{3,2}^{g+1} + 64\cdot 15 S_{4,1}^{g+1})+\right.$$
$$\left. 24(14  S_{3,2}^{g+1} +
 112\cdot 15 S_{4,1}^{g+1} +2^9\cdot 155 S_{5, 1/2}^{g+1})\right)=$$
 $${2^{-2g+7}}(S_{3,2}^{g+1}+90S_{4,1}^{g+1}+3720 S_{5,1/2}^{g+1}).$$
Rescaling from $g+1$ to $g$, and gathering all $S_{3,2}^g$ on one side, we get
\begin{prop}
For any $g\geq 5$ and for any $\pi\in\T_g$ we have
\begin{equation}\label{Sexpress}
 (2^{2g-6}-1)S_{3,2}^g(\pi)= 90S_{4,1}^g(\pi)+3720 S_{5,1/2}^g(\pi)
\end{equation}
(notice that the last term $S_{5,1/2}^g$ is only known by \cite{RSM} to be a modular form on $\T_g$, so the above identity does not make sense over all of $\H_g$).
\end{prop}
Substituting this result in formula (\ref{Xi}) for $\Xi^{(5)}(\pi)$ for $\pi\in\T_5$ to express $S_{5,1/2}^5$ in terms of $S_{4,1}^5$ and $S_{3,2}^5$, and then using lemma \ref{iglm} to express those, we eventually express $\Xi^{(g)}$ as a linear combination of $S_{0,16}^5$ and $S_{1,8}^5$, and furthermore as a linear combination of $(f_4^{(g)})^2$ and $f_8^{(g)}$, cf.~\cite{Ch}. We thus obtain
\begin{thm}\label{XiF}
For any $\pi\in\T_5$ we have
$$
 \Xi^{(5)}(\pi)=\frac{-51}{217}((f_4^{(5)})^2-f_8^{(5)})(\pi)
$$
\end{thm}
\begin{cor}
For $g= 5$ the cosmological constant $\Xi^{(5)}$ vanishes identically on $\T_5$ if and only if $(f_4^{(5)})^2-f_8^{(5)}$ vanishes identically on $\T_5$.
\end{cor}
It was conjectured in \cite{BK,BK1,M86,159} that $(f_4^{(g)})^2-f_8^{(g)}$ vanishes identically on $\T_g$ in any genus (physically $f_4$ is interpreted as the appropriate measure for the $E_8$ theory, $f_4^2$ --- for $E_8\times E_8$, and $f_8$ --- for $SO(32)$). For $g\le 3$ the identical vanishing of this modular form  on $\T_g=\H_g$ is a consequence of Riemann's bilinear addition theorem. For $g=4$ this form is equal to the Schottky equation defining $\T_4\subset\H_4$, cf\cite{Ch} . In cf\cite{Po} it was shown that this form vanishes along the hyperelliptic locus for any $g$. We will now show that in fact this form does {\it not} vanish identically on $\T_5$, from which it follows that it does {\it not} vanish identically on $\T_g$ for any $g\ge 5$.

\section{The generic non-vanishing of $(f_4^{(5)})^2-f_8^{(5)}$ along $\T_5$}
Writing out the formulas for $S_{0,16}^g$ and $S_{1,8}^g$ explicitly, we have from (\ref{fS})
$$
 F_g(\tau):=-2^{2g}\left((f_4^{(g)})^2(\tau)-f_8^{(g)}(\tau)\right)=2^g\sum\limits_{m\in\F^{2g}}
 \t_m^{16}(\tau)-\left(\sum\limits_{m\in\F^{2g}} \t_m^8(\tau)\right)^2
$$

We now prove that the cosmological constant $\Xi^{(5)}(\tau)$ for the ansatz above does not vanish identically for $\tau\in\T_5$. By theorem \ref{XiF} it is equivalent to proving the following for $g=5$.

\begin{thm}\label{nonvanish}
The modular form $F_g$ does not vanish identically on $\T_g$.
\end{thm}

\begin{proof}
We will prove this by showing that $F_5$ does not generically vanish in a neighborhood of the boundary divisor $\delta_0\subset\overline{\M_5}$ (or does not vanish along $\delta_0$ to second order), and then deducing that $F_g$ does not vanish identically for any $g\ge 5$.

Indeed, consider a family of genus 5 Riemann surfaces $C_t$ degenerating to $(C,p,q)/p\sim q\in \M_{4,2}=\delta_0\subset\overline{\M_5}$. In such a family, by \cite[Corollary 3.8, p.~53]{faybook} that the period matrix $\tau_t$ of $C_t$ has the form
$$
 \left(\begin{matrix} w(t) &z(t)'\\ z(t)& \tau(t)\end{matrix} \right)=\left(\begin{matrix} \ln t+ c_1+ c_2t&A_{pq}'+\frac14 tv_{pq}'\\ A_{pq}+\frac14 tv_{pq}& \tau+ t \sigma\end{matrix} \right),
$$
where $\tau$ is the period matrix of $C$, $A:C\to Jac(C)$ is the Abel-Jacobi map, $v:C\to\PP^{g-1}$ is the associated canonical curve (Gauss map for $A$), $A_{pq}:=A(p)-A(q)$,
$v_{pq}:=v(p)-v(q)$, ' denotes the transpose, and we have set
$\sigma_{ij}:=\frac{1}{4} (v_i(p)- v_i(q))(v_j(p)- v_j(q))$

Recall now the Fourier-Jacobi expansion of the theta functions near the boundary, see\cite{vG}:
$$
\tt{0\ \e}{\de_1\ \de}\left(\begin{matrix}\tau_{11}&z^t\\ z& \tau\end{matrix}\right)=
 \tt\e\de(\tau,0)+2e^{\pi i\de_1}q^4\tt\e\de(\tau,z)+O(q^{16})
$$
and
$$
\tt{1\ \e}{\de_1\ \de}\left(\begin{matrix}\tau_{11}& z^t\\ z& \tau\end{matrix}\right)=2e^{\pi i\de_1/2}q\tt\e\de(\tau,z/2)+O(q^9)
$$
where as usual we let $q:=\exp(\pi i \tau_{11}/4)$; to relate this to Fay's notation above, we have $t=q^8$.

Let us now compute the first terms of the Taylor expansion of $F_g$ as $t\to 0$ (i.e. near the boundary, as $\tau_{11}\to i\infty$). By inspection we see that the two lowest order terms are $O(1)$ and $O(q^8)=O(t)$ respectively, so we compute them using
$$
 \sum\limits_{\e,\de\in\F^{g-1}}\sum\limits_{\de_1\in\F} \theta^N\tch{0\ \e}{\de_1\ \de}\left(\begin{matrix}\tau_{11}&z^t\\ z& \tau\end{matrix}\right)
$$
$$=2\sum\limits_{\e,\de\in\F^{g-1}}
 \theta^N\tch\e\de(\tau,0)+2\binom{N}{2}(2q^4)^2\theta^{N-2}\tch\e\de(\tau,0)\theta^2
 \tch\e\de(\tau,z)+o(q^8)
$$
and
$$
 \sum\limits_{\e,\de\in\F^{g-1}}\sum\limits_{\de_1\in\F} \theta^N\tch{1\ \e}{\de_1\ \de}\left(\begin{matrix}\tau_{11}&z^t\\ z& \tau\end{matrix}\right)=2^Ne^{N\pi i\de_1/2}q^N\theta^N\tch\e\de(\tau,z/2)+O(q^N).
$$
Adding these together, we get the Fourier-Jacobi expansions
$$
 \sum\limits_{\e,\de\in\F^g}\theta^{16}\tch\e\de(\tau_t)\left(\begin{matrix}\tau_{11}&z^t\\ z& \tau\end{matrix}\right)
$$
$$
 =2\sum\limits_{\e,\de\in\F^{g-1}{\rm\ even}}
 \theta^{16}\tch\e\de(\tau,0)+2\binom{16}{2}(2q^4)^2\theta^{14}\tch\e\de(\tau,0)\theta^2
 \tch\e\de(\tau,z)+o(q^8)
$$
with no contribution from the case of $\e_1=1$, while
$$
 \sum\limits_{\e,\de\in\F^g}\theta^{8}\tch\e\de\left(\begin{matrix}\tau_{11}&z^t\\ z& \tau\end{matrix}\right)=2\sum\limits_{\e,\de\in\F^{g-1}}
 \theta^{8}\tch\e\de(\tau,0)+2\binom{8}{2}(2q^4)^2\theta^{6}\tch\e\de(\tau,0)\theta^2
 \tch\e\de(\tau,z)
$$
$$+\sum\limits_{\e,\de\in\F^{g-1}}2^8e^{8\pi i\de_1/2}q^8\theta^8\tch\e\de(\tau,z/2)+o(q^8).
$$
Combining these, and substituting $\tau+t\sigma$ from Fay's formula for $\tau$ in the Fourier-Jacobi expansion above we get for the lowest order terms of $F_g(\tau_t)$ the expression
$$
 F_5(\tau_t)= 4F_{4}(\tau+t\sigma)+2^5\cdot960t\sum\limits_{\e,\de\in\F^{g-1}}
 \theta^{14}\tch\e\de(\tau+t\sigma,0)\theta^2\tch\e\de(\tau+t\sigma,z)
$$
$$
 -2t\!\!\!\sum\limits_{\alpha,\beta\in\F^{g-1}}\!\!\!
 \theta^8\tch\alpha\beta(\tau+t\sigma,0) \left(
 448\!\!\!\sum\limits_{\e,\de\in\F^{g-1}}\!\!\!
 \theta^6\tch\e\de(\tau+t\sigma,0)\theta^2\tch\e\de(\tau+t\sigma,z)\right.
$$
$$
 \left.
 +512\sum\limits_{\e,\de\in\F^{g-1}}
 \theta^8\tch\e\de(\tau+t\sigma,z/2)\right)+o(t)
$$
where we denoted $z=A_{pq}+\frac14 tv_{pq}$, and recall that $t=q^8$.

We are now ready to prove that $F_5$ does not vanish identically on $\M_5$. Indeed, if it vanished identically, it would vanish along $\delta_0$ to any order in the expansion. However, notice that the linear in $t$ term of the expansion of $F_5$ has two summands --- one obtained by setting $t=0$ in the term multiplying $t$ above, and one obtained by expanding $4F_4(\tau+t\sigma)$, i.e.~we have
$$
 F_5(\tau_t)= 4F_{4}( \tau) +t ( G(\tau, \sigma)+  f_1(\tau, (A(p)-A(q)))+o(t)
$$
where the term $G(\tau, \sigma)$ is given by
$$G(\tau, \sigma)=\sum_{i\leq j}\frac{\partial F_4}{\partial\tau_{ij}}\sigma_{ij}(p, q)$$
and $f_1$ denotes the term multiplying $t$ above, evaluated at $t=0$

Our goal now is to show that $G+f_1$ does not vanish identically.
One potential argument to show this would be to view $G$ and $f_1$ as functions of $p$ and $q$ on the universal cover of $C\times C$, as $p$ and $q$ move. If one could show that they are sections of different line bundles, then their sum would vanish identically if and only if they both vanished identically. However, to do this rigorously requires a very careful analysis of the above degeneration on the universal cover, and the details would be very messy.

Instead, we first use Riemann relations to simplify $f_1$, by expressing $\theta^8\tch\e\de(\tau,z/2)$ above as a linear combination of $\theta^2\tch\e\de(\tau,z)$, and then use the beautiful genus 4 geometry to consider the behavior of the $G+f_1$ term as $p$ and $q$ approach each other.
\begin{lm}
We have the following identity:
$$
 \sum\limits_{\e,\de\in\F^{g-1}}\theta^8\tch\e\de(\tau,z/2)
=\sum\limits_{\e,\de\in\F^{g-1}}\theta^6\tch\e\de(\tau,0)\theta^2\tch\e\de(\tau,z).
$$
\end{lm}
\begin{proof}
A special case of Riemann relations (\ref{riemannrelation}) is
$$
\theta^4\tch\e\de(\tau,z/2)=2^{1-g}\sum\limits_{\alpha,\beta \in\F^{g-1}} (-1)^{\alpha\cdot\de+\beta\cdot\e}\theta^3\tch\alpha\beta(\tau,0)\tt\alpha\beta(\tau,z).
$$
We will now use this identity twice to get an expression for $\theta^8\tch\e\de(\tau,z/2)$ as a double sum, and then sum over all $\e,\de\in\F^{g-1}$ (note that we include the odd ones), to get
$$
 \sum\limits_{\e,\de\in\F^{g-1}}\theta^8\tch\e\de(\tau,z/2)
 =2^{2-2g}\sum\limits_{\epsilon,\delta,\alpha,\beta,\sigma,\mu\in\F^{g-1}}
$$
$$
 (-1)^{(\alpha+\sigma)\cdot\de+(\beta+\mu)\cdot\e}\theta^3\tch\alpha\beta(\tau,0)\tt\alpha\beta(\tau,z)
 \theta^3\tch\sigma\mu(\tau,0)\tt\sigma\mu(\tau,z).
$$
Notice that the only dependence on $\e,\de$ in the sum on the right is in the sign. Recalling that in general for any $B\in\F^g$
$$
 \sum\limits_{A\in\F^g}(-1)^{A\cdot B}=2^g\de_{B,0}
$$
(where $\de_{B,0}$ is the Kronecker symbol), we see that for fixed $\alpha,\beta,\sigma,\mu$ the sum over $\e,\de$ on the right-hand-side of the formula above is non-zero if and only if $\alpha+\sigma=\beta+\mu=0$ (i.e. iff $\alpha=\sigma$ and $\beta=\mu$). The factor of $2^{2g-2}$ for summing over $\e,\de$ in this case cancels out the $2^{2-2g}$ from Riemann relations, so that we finally obtain
$$
 \sum\limits_{\e,\de\in\F^{g-1}}\theta^8\tch\e\de(\tau,z/2)
=\sum\limits_{\alpha,\beta\in\F^{g-1}}\theta^6\tch\alpha\beta(\tau,0)\theta^2\tch\alpha\beta(\tau,z)
$$
as claimed.
\end{proof}

Using the lemma, the $tf_1$ term becomes (note that $448+512=960$)
$$
 1920t\left(2^{g-1}\sum\limits_{\e,\de\in\F^{g-1}}
 \theta^{14}\tch\e\de(\tau,0)\theta^2\tch\e\de(\tau,z)\right.
$$
$$\left.-
 \sum\limits_{\alpha,\beta\in\F^{g-1}}
 \theta^8\tch\alpha\beta(\tau,0)
 \sum\limits_{\e,\de\in\F^{g-1}}
 \theta^6\tch\e\de(\tau,0)\theta^2\tch\e\de(\tau,z)\right)
$$
If we now denote, in the spirit of \cite{vgvdg}, the variables
$$
 X_{\e,\de}:=\theta^2\tch\e\de(\tau,0)
$$
and express $F_{g-1}$ as a function of $X$'s, note that the above expression is
$$
 f_1(\tau,z)=v_{F_{g-1}}:=\sum\limits_{\e,\de\in\F^{g-1}}240 \frac{\partial F_{g-1}}{\partial X_{\e,\de}}\ \theta^2\tch\e\de(z)
$$
(note that when differentiating $X_{\e,\de}^8=\theta^{16}\tch\e\de(\tau,0)$, one picks up a factor of 8, and while differentiating $(\sum X_{\e,\de}^4)^2$, one also picks up the factor of $2\cdot 4=8$ for each term).

\smallskip
We will now deal with the term $G$, and investigate $G+f_1$ as the points $p$ and $q$ approach each other. Indeed, suppose that in a local coordinate $z$ on a chart on a Riemann surface  $z(p)=z(q)+u$, and consider the expansion of this linear in $t$ term as a function of $u$ near $u=0$. We have
$$
 G(\tau, \sigma)=u^2 \frac{1}{4} \sum_{i\leq j}\frac{\partial F_4}{\partial\tau_{ij}}v_i'(q) v_j '(q)+ u^3  \frac{1}{2} \sum_{i\leq j}\frac{\partial F_4}{\partial\tau_{ij}}v_i''(q) v_j '(q)+ O(u^4).
$$
Since $f_1(\tau,z)$ is an even function of $z$, from the above formula we have
$$
 f_1(\tau, A(z(q)+u)-A(z(q))) =u^2 (2\pi i) 30  \sum_{i\leq j}\frac{\partial F_4}{\partial\tau_{ij}}v_i(q) v_j (q)+O(u^4).
$$

Since $F_4$ is the Schottky form, which gives the defining equation for $\M_4\subset\A_4$, the partial derivatives  $\frac{\partial F_4}{\partial\tau_{ij}}$ give the tangent directions to $\M_4$. By the results of Andreotti and Mayer \cite{anma}, they are proportional to the coefficients of the quadric $Q$ containing the canonical curve.

Thus for the  $u^2$ coefficient
$$
 \frac{1}{4} \sum_{i\leq j}\frac{\partial F_4}{\partial\tau_{ij}}v_i'(q) v_j '(q)+ (2\pi i) 30  \sum_{i\leq j}\frac{\partial F_4}{\partial\tau_{ij}}v_i(q) v_j (q)=0
$$
The second summand vanishes identically, and thus we would need to have
$$
 \sum_{i\leq j}\frac{\partial F_4}{\partial\tau_{ij}}v_i'(q) v_j '(q)\equiv 0
$$
which is to say that $Q$ would contain the Gauss map image of the canonical curve, i.e.~$^t v'(q)Qv'(q)$  would be zero for all $q$.

Now consider the $u^3$ coefficient of the expansion; since $f_1$ is an even function of $z$, it is equal to
$$
  \frac{1}{2} \sum_{i\leq j}\frac{\partial F_4}{\partial\tau_{ij}}v_i''(q) v_j '(q)=\frac12 ^t v''(q)Qv'(q)
$$
Since the quadric contains the curve, along which $u$ is the parameter $u$,  we would also have
$$
 ^t v(q+u)Qv(q+u)=0\qquad \forall q\in \mathcal C, \forall u\in T_q C
$$
Expanding this in $u$ would yield
$$
 0=^t v(q)Qv(q)+ 2u  ^t v'(q)Qv(q) +u^2 ( 2  ^t v''(q)Qv(q) +     ^t v'(q)Qv'(q) )+
$$
$$
 2u^3 ( ^t v''(q)Qv'(q)+^t v'''(q)Qv(q)) + O(u^4)
$$
and thus, since the form was assumed to vanish identically,
$$
 0=v(q)Qv(q)=v(q)Qv'(q)=v(q)Qv''(q)=v(q)Qv'''(q)
$$
However, $Qv(q)$ is generically non-zero, since the canonical curve is non-degenerate in ${\mathbb P}^3$. Hence the matrix
$$(v(q), v'(q),v''(q),Qv'''(q))$$
must have non-maximal rank (i.e.~rank at most 3). However, this matrix is just the Wronskian (here we are using $g=4$), and thus $q$ must be a Weierstrass point, while we assumed that the above must hold for all points
$q\in \mathcal C$. We have thus arrived at a contradiction, and thus $F_5$ does not vanish identically.

To show now that $F_g$ does not vanish identically on $\T_g$ for any $g>5$, note that in the above expansion of $F_g$ along $\delta_0$ the constant term is $4F_{g-1}$, so if $F_{g-1}$ does not vanish identically on $\T_{g-1}$, then $F_g$ cannot vanish identically on $\T_g$.
\end{proof}

Note that for $g=5$ the slope of the Brill-Noether divisor is $6+\frac{12}{g+1}=8$. Since the slope conjecture is valid for $\overline{\M_5}$ (see cf.~\cite{FP} for more discussion), this is the unique effective divisor on $\overline{\M_5}$ of slope 8. However, the slope of the divisor of $F_5$ is equal to 8 as well.
\begin{cor}
The zero locus of the modular form $F_5$ on $\M_5$ is the Brill-Noether divisor or, equivalently, the locus of trigonal curves.
\end{cor}
It is interesting to ask whether $F_g$ might vanish on the locus of trigonal curves for any genus (note that it is known that $F_g$ vanishes on the locus of hyperelliptic curves).
\begin{rem}
We recall that there exists the Prym map $p:{\mathcal R}_5\to\A_4$ from the moduli space of curves of genus 5 with a choice of a point of order two. Let us also denote by $q:{\mathcal R}_5\to\M_5$ the forgetful finite map. From Recillas construction it then follows, cf.~\cite{BL}, that the trigonal locus in $\M_5$ is in fact equal to $q\circ p^{-1}(\M_4)$, and it would be interesting to understand the vanishing of $F_5$ in these terms, by exploring the Schottky-Jung identities more explicitly.
\end{rem}

\begin{rem}
In a forthcoming paper \cite{OPSMY} the space of cusp forms $[\Gamma_4(1,2),8)]_0$, in which $\Xi^{(4)}[0]$ lies, is studied, and it is shown that the dimension of this space is equal to 2. By looking at the basis for this space it follows that there exists a unique cusp form with prescribed factorization properties, and thus it follows that $\Xi^{(4)}[0]$ is the unique form in $[\Gamma_4(1,2),8]_0$ satisfying the factorization constraints.
\end{rem}

\begin{rem}
In \cite{FD} 7 linearly independent modular forms in $[\Gamma_4(1,2),8]$ are constructed, which are polynomials in squares of theta constants, and $\Xi^{(4)}[0]$ is expressed as their linear combination. If this expression for $\Xi^{(4)}[0]$ is in fact divisible by $\theta_0^2(\tau,0)$, then it would follow that in the expression for the genus 4 two-point function
$$
 \sum_m \frac{ \Xi^{(4)}[m](\tau)} {\t_m^2(\tau,0)}\t_m^2(\tau,z)
$$
all ratios are polynomials in squares of theta constants. It is tempting then to conjecture that this two-point function is a constant multiple of $v_{F_4}$, and thus the vanishing of the genus 4 two-point function would be equivalent (see \cite{GSM}) to the function $v_{F_4}(\tau)$ lying in $\Gamma_{00}$ for any $\tau\in\T_4$, which by the above discussion is {\it not} the case. This question merits further investigation.
\end{rem}
In contrast to genera up to 4, it turns out that in genus 5 the modular forms $\Xi^{(5)}[m]$ are not the unique forms on $\T_5$ satisfying the factorization constraints.
\begin{prop}\label{new}
For any constant $c$ (independent of $m$) the expressions $$\Xi'^{(5)}[m]:=\Xi^{(5)}[m]+c(f_4^2-f_8)$$ are modular form of weight 8 on $\T_5$ with respect to the subgroup of $\operatorname{Sp}(5,\Z)$ fixing $m$, permuted among themselves among themselves under the action of $\operatorname{Sp}(5,\Z)$ satisfying factorization constraints. Moreover, for $c=\frac{38192}{17}$ the cosmological constant $\sum_m\Xi'^{(5)}[m]$ vanishes identically on $\T_5$.
\end{prop}
\begin{proof}
Note first that $f_4^2-f_8$ is a modular form of weight 8 with respect to all of $\operatorname{Sp}(5,\Z)$, and thus $\Xi'^{(5)}[m]$ are modular and permuted by the group action as claimed. To determine the factorization of $\Xi'^{(5)}[m]$ on $\T_i\times\T_{5-i}$, we compute in general for $\tau_1\in\T_i,\tau_2\in\T_{g-i}$
$$
 f_4^2\left(\begin{matrix}\tau_1&0\\ 0&\tau_2\end{matrix}\right)
 -f_8\left(\begin{matrix}\tau_1&0\\ 0&\tau_2\end{matrix}\right)
 =f_4^2(\tau_1)f_4^2(\tau_2)-f_8(\tau_1)f_8(\tau_2).
$$
If $i\le 4$ (which is always the case for $g=5$), so that $F_i=0$, the above expression becomes equal to
$$
 f_8(\tau_1)(f_4^2(\tau_2)-f_8(\tau_2));
$$
if now $g-i\le 4$ (which is also always the case for $g=5$), this vanishes, and thus the factorization of $\Xi'$ is the same as the factorization of $\Xi$.

Finally we note that by definition the cosmological constant is
$$
 \sum\limits_m\Xi'^{(5)}[m]= \sum\limits_m\Xi^{(5)}[m]+2^4(2^5+1)(f_4^2-f_8)
 \left(-\frac{51}{217}+528 \right)(f_4^2-f_8)
$$
where we used the computation of the cosmological constant for $\Xi^{(5)}$ from theorem \ref{XiF}, and recall that the number of even characteristics is $2^{g-1}(2^g+1)$. Thus for $c=\frac{38192}{17}$ the cosmological constant for $\Xi'^{(5)}$ vanishes identically.
\end{proof}

\section{Appendix: possible generalizations to higher genus}
The results of the previous section lead us to observe that a way to have modularity in the ansatz for the chiral superstring measure proposed in \cite{GR} for all $g$ is to prove a generalized version of Schottky-Jung relations, involving roots of degree $2^{k-4}$ for all $5\le k\le g$. These generalized Schottky-Jung relations should be relations induced by the Riemann relations in genus $g-k+4$. For example, cf.~\cite{AC}, a Riemann relation in genus 4
$$r_1\pm r_2\pm r_3\pm r_4=0,$$
where each $r_i$ is of the form
$$
 r_i=\theta_{m_1}\theta_{m_2}\theta_{m_3}\theta_{m_4},
$$
induces a Schottky relation in genus $6$ of the form
$$
 S_1\pm S_2\pm S_3\pm S_4=0
$$
where each $S_i$ is a fourth root of a monomial of degree 16 in the theta constants of a Jacobian of a genus 6 curve, with the set of characteristics satisfying some obvious conditions.

If this is the case, and such generalized Schottky-Jung relations hold, as immediate consequence of Riemann's formula, we have the following relations
\begin{prop}\label{Srelation}
If the generalized Schottky-Jung relations hold, then for any $g\geq k+2$, for any $\tau\in\T_g$, we have
$$
 (2^{2g-2k}-1)S_k^g=6(2^{k+1}-1) S_{k+1}^g+ 8(2^{k+2}-1)(2^{k+1}-1) S_{k+2}^g
$$
(this is a generalization of (\ref{Sexpress}).
\end{prop}
If this is the case, then by eliminating $S_k^g$ starting from the highest one, $S_g^g$, the cosmological constant $\Xi^{(g)}$ given by (\ref{Xi}) on $\T_g$ can be written as a linear combination of $S_{0,16}^g$ and $S_{1,8}^g$, or as a linear combination of $(f_4^{(g)})^2$ and $f_8^{(g)}$.
Thus it makes sense to ask if $\Xi^{(g)}$ is proportional to the restriction of $(f_4^{(g)})^2-f_8^{(g)}$ also when $g>5$. Eberhard Freitag confirmed this, using a computer, for $g<1000$. We now give a (geometric, rather than combinatorial) proof of this.
\begin{prop}
If proposition \ref{Srelation} holds, then we have for any $\tau\in\T_g$
$$\Xi^{(g)}(\tau)={\rm const}(f_4^2(\tau)-f_8(\tau)).$$
\end{prop}
\begin{proof}
From the above discussion we must have
$$\Xi^{(g)}= a_g(f_4^{(g)})^2- b_g f_8^{(g)}$$
for some constants $a_g$ and $b_g$. Let us introduce the Siegel $\Phi$-operator: for any $f:\H_g\to \C$ we let
$$
 \Phi(f)(\tau_1):=\lim_{ \lambda \longrightarrow + \infty}f \pmatrix \tau_1 &0\\ 0&i\lambda\endpmatrix
$$
for all $\tau_1\in\H_{g-1}$. This operator has relevance in the theory of modular forms, cf.~\cite{I,Fr} for details. Applying the Siegel $\Phi$ operator to forms defined on $\T_g$ we get forms defined on $\T_{g-1}$. It is well known and easy to show using the expression in terms of $S_0$ and $S_1$ that
$$
 \Phi(f_4^{(g)})=f_4^{(g-1)}\quad{\rm and } \quad \Phi(f_8^{(g)})=f_8^{(g-1)}
$$
An easy computation, cf.~\cite{RSM}, then gives $\Phi(\Xi^{(g)})=0$, and moreover
$$
 \Phi^{g-4}(a_g(f_4^{(g)})^2- b_g f_8^{(g)})=a_g(f_4^{(4)})^2- b_g f_8^{(4)}.
$$
Thus we have
$$
 0= \Phi^{g-4}(\Xi^{(g)})=a_g(f_4^{(4)})^2- b_g f_8^{(4)}.
$$
This implies that $a_g=b_g$, since $(f_4^{(4)})^2-f_8^{(4)}$ is the defining equation for $\T_4\subset\H_4$. Hence setting $c_g:=a_g=b_g$, we have the desired equality
$$
 \Xi^{(g)}= c_g((f_4^{(g)})^2- f_8^{(g)})
$$
for all $g$.
\end{proof}
\begin{cor}
If the generalized Schottky-Jung identities hold (i.e. if proposition \ref{Srelation} holds) $\Xi^{(g)}$ vanishes identically on $\T_g$ if and only if $g\le 4$.
\end{cor}
\begin{rem}
It is tempting to try to construct, similarly to proposition \ref{new}, a corrected ansatz for arbitrary genus that would satisfy (assuming the above holds) the factorization constraints and give a vanishing cosmological constant. However, already in genus 6 it is not clear how to proceed. It is natural to try to add a multiple of $f_4^2-f_8$ to $\Xi^{(6)}$, but similarly to the proof of proposition \ref{new} we see that in this case on $\T_1\times\T_5$ the term  $f_8^{(1)}F_5$ is added to the factorization, which is proportional to $f_8^{(1)}(\Xi'^{(5)}[m]-\Xi^{(5)}[m])$, but not proportional to the necessary $\Xi^{(1)}[m](\Xi'^{(5)}[m]-\Xi^{(5)}[m])$.
\end{rem}

\end{document}